# An algorithm for the evaluation of the Gamma function and ramifications

D. KARAYANNAKIS
T.E.I. of Crete
Dept. of Science / Applied Math. Sector
dkar@stef.teiher.gr

Abstract
A closed type formula is provided for Γ(x), when x is rational in (0,1), based on a recursive algorithm which, along with its equivalent formulation as an infinite product, furnishes easy to obtain estimates and leads to an additional insight on various properties concerning Γ(x), 1/Γ(x), lnΓ(x), B(x,y), ψ(x) and ψ'(x).The method allows the production of estimates of non-classical values of Γ(x) and a variety of pertinent facts without resorting to any kind of specialized programming tools beyond a scientific pocket calculator.Formulae some of them seemingly new, involving infinite products and inequalities, are presented as well.

## 1. INTRODUCTION

The theory of the Gamma function, Γ(x), was historically developed in connection with the problem of extending the value of n! for positive integers n to arbitrary real numbers.The most fundamental functional equation, that can be derived almost immediately from its definition and which is basic for the development of the theory of the Gamma function is:

$$\Gamma(x+1) = x\,\Gamma(x). \qquad (1.1)$$

It represents a generalization of the identity n!=n(n-1)! for non integral values of n. Note that Eq.(1.1) allows us to extend the original classical definition to include negative real numbers, provided we avoid negative integral values and zero. In the rest of the paper we will always assume that knowledge of the values of the Gamma function on the interval (0, 1] results to knowing the values of Γ(x) on $\Re \setminus \mathbb{Z}_0^-$.

Legendre's doubling formula states that ([4]):

$$\Gamma(2x) = \frac{2^{2x-1}}{\pi^{1/2}} \Gamma(x+1/2)\Gamma(x). \qquad (1.2)$$

When (1.2) is employed for x=1/2 it gives immediately Γ(1/2)=π**1/2**, while for x=n can lead to the evaluation of Γ(n+1/2) in closed form (cf. [6]); but these are basically the most this formula can do to that direction.

Eq. (1.3) is a special case of Gauss' multiplication formula ([4]):

$$\Gamma(nx) = \frac{n^{(nx-1/2)}}{(2\pi)^{(n-1)/2}} \Gamma(x)\,\Gamma(x+1/n)\ldots\Gamma(x+(n-1)/n). \qquad (1.3)$$

Eq. (1.3) in most of the cases leads to identities that at first glance look rather exotic (cf. [6]) and its use is subjected to the same, more or less, limitations described above for Eq. (1.3); nevertheless it will be instrumental for our investigations since it will provide the constant that multiplies the infinite product representation of Γ(1/p), for any integer p, p≥3 (Lemma 3.1).

## 2. A CLOSED TYPE FORMULA FOR $\Gamma(x+b)/\Gamma(x)$, $x>0$, $0 \leq b \leq 1$

It is known that a hypergeometric series is a series of the form:

$$F(\alpha, \beta, \gamma; z) = 1 + \frac{\alpha \cdot \beta}{\gamma \cdot 1} z + \frac{\alpha \cdot (\alpha+1) \cdot \beta \cdot (\beta+1)}{\gamma \cdot (\gamma+1) \cdot 1 \cdot 2} z^2 + \ldots \tag{2.1}$$

A hypergeometric series terminates if $\alpha$ or $\beta$ is a negative integer or zero. If we exclude the case $\gamma = -n$, $n=0,1,2,\ldots$ the series (2.1) converges in the unit circle $|z|<1$. The sum-function, the hypergeometric function F, (where for simplicitly we prefer this symbolism instead of the usual 2F1) has a branch point at z=1. Still, whenever $Re(\alpha+\beta-\gamma)<0$, the series converges absolutely throughout the entire unit circle ([10]).

From now on we will consider only the case of real parameters $\alpha$, $\beta$, $\gamma$ and real z= t and we will be concerned mainly with the special case t=1. We have the following formula known as Gauss Hypergeometric (or Summation) Theorem ([10]):

$$F(\alpha, \beta, \gamma; 1) = \frac{\Gamma(\gamma)\Gamma(\gamma-\alpha-\beta)}{\Gamma(\gamma-\alpha)\Gamma(\gamma-\beta)} \tag{2.2}$$

under the condition $\gamma>\alpha+\beta$; (2.2) is valid also for complex variables whose real parts obey the analogous condition.

It is also known ([2]) that the (real) hypergeometric series is one of the solutions of the differential equation:

$$t(1-t)\, y''(t) + [\gamma - (\alpha+\beta+1)t]\, y'(t) - \alpha\beta y(t) = 0. \tag{2.3}$$

We recall that for the function psi, i.e. the logarithmic derivative of $\Gamma(x)$, we have the formula ([6]):

$$\psi(u) - \psi(t) = \sum_{k=0}^{\infty} \left( \frac{1}{t+k} - \frac{1}{u+k} \right) \tag{2.4}$$

whenever u and t are nonnegative integers or zero.

We will first present a "formal" (heuristic) proof of the following:

LEMMA 2.1. $F(\alpha, b, 1; t) = \exp\left[ \sum_{n=1}^{\infty} g_n(\alpha,b) t^n \right]$ for $\alpha+b<1$, $0 \leq t \leq 1$, where $g_1(\alpha,b)=\alpha b$, $g_2(\alpha,b)=\alpha b(1+\alpha+b-\alpha b)/4$, and the $g_n$'s (for $n \geq 2$) are defined through the recursive algorithm:

$$(n+1)^2 g_{n+1} = n(n+\alpha+b-g_1)g_n + \sum_{k=0}^{n-2} (k+1)g_{k+1}\left[(n-k-1)g_{n-k-1} - (n-k)g_{n-k}\right]. \tag{2.5}$$

PROOF. The differential equation (2.3) for $\gamma=1$ leads, after the substitution:

$$F(\alpha, b, 1; t) = C_1 \exp g(t)$$

to the equation:

$$t(1-t)[(g')^2 + g''] + [1-(\alpha+b+1)t]g' - \alpha b = 0. \tag{2.6}$$

Solving (2.6) formally by the power-series method ([8]) we obtain $g(t)=g_0 + \sum_{n=1}^{\infty} g_n(\alpha,b)t^n$, where $g_n(\alpha,b)$ can be retrieved via (2.5) for $n \geq 2$. But, since:

$$F(\alpha, b, 1; 0)=1$$

we can conclude immediately that $k \exp g_0 = 1$, which leads to the announced result.

REMARK 2.1. From now on we will consider $1-\alpha-b=x>0$ and then in Lemma 2.1, by some abuse of notation, the functions $g_n$ will be considered functions of x and b, and, in particular, $g_1=b(1-b-x)$ and $g_2 = b(1-b-x)\dfrac{(b-1)x + b^2 - b + 2}{4}$. Observe also that the coefficient of $g_n$ in (1.5) takes now the form $n[n+1-x-b(1-b-x)]$.

Combining now Lemma 2.1 and (2.2) we can conclude immediately:

PROPOSITION 2.1. $\Gamma(x+b) = f(x,b)\dfrac{\Gamma(x)}{\Gamma(1-b)}$, for $x > 0$, $0 \leq b < 1$, where

$$f(x,b)=\exp\left[-\sum_{n=1}^{\infty}g_n(x,b)\right].$$

COROLLARY 2.1. $f(x,b) \sin(\pi b) = f(b,x) \sin(\pi x)$, $0 < x, b < 1$.

PROOF. Using the formula of Proposition 2.1 we automatically see that for $0<x, b<1$
$\Gamma(x+b)=f(x,b)\dfrac{\Gamma(x)}{\Gamma(1-b)} = \Gamma(b+x) = f(b,x)\dfrac{\Gamma(b)}{\Gamma(1-x)} \Rightarrow f(x,b)\, \Gamma(x)\, \Gamma(1-x) = f(b,x)\, \Gamma(b)\, \Gamma(1-b)$.
But since $\Gamma(z)\, \Gamma(1-z) = \dfrac{\pi}{\sin(\pi z)}$ (cf. [4]) for $0<z<1$, we arrive immediately at the announced result.

We are now ready to prove a theorem, which will allow us to recapture known results, as well as to produce infinite product representations of values of $\Gamma(x)$ with a satisfactory precision by the mere use of a pocket calculator (see Appls 9 and 10 in Ch. 4 as well as the results in [7]).

For reasons that should be already apparent through the functional relation of Proposition 2.1, and later on through Theorem 3.1 which is the main result of Ch. 3, we will name the bivariable function $f(x,b)$ the <u>joint factor</u> of the $\Gamma(x+b)$ function.

THEOREM 2.1. The joint factor has the following infinite-product representation:

$$f(x,b)=\prod_{k=1}^{\infty}\frac{k(x+k-1)}{(k-b)(x+k+b-1)}. \qquad (2.7)$$

PROOF. Fix temporarily b in (0,1). Then differentiating for x the functional equation of Proposition 2.1, we obtain: $\Gamma(x+b)\,\psi(x+b) = f(x,b)\,\Gamma(x) + f(x,b)\,\Gamma(x)\,\psi(x)$ and thus

$f(x,b)\Gamma(x)\psi(x+b) = -f(x,b)\left(\sum_{n=1}^{\infty}g'_n(x,b)\right)\Gamma(x) + f(x,b)\Gamma(x)\psi(x)$. We conclude that

$\sum_{n=1}^{\infty}g'_n(x,b)= \psi(x)-\psi(x+b)$. Appealing now to (2.4) for $t=x+b$ and $x=u$, we obtain that:

$\sum_{n=1}^{\infty}g'_n(x,b) = \sum_{k=0}^{\infty}\left(\dfrac{1}{x+b+k} - \dfrac{1}{x+k}\right)$. Integrating now for x we obtain directly that: $\sum_{n=1}^{\infty}g_n(x,b)=C_2+\sum_{k=0}^{\infty}\ln\dfrac{x+b+k}{x+k}$. Observing now that $g_1(x,1-b)=g_2(x,1-b)=0$, the recursive formula (2.5) leads to $g_n(x,1-b)=0$ for all n. We conclude that

$C_2=-\sum_{k=0}^{\infty}\ln\dfrac{1+k}{1-b+k}$ and thus $\sum_{n=1}^{\infty}g_n(x,b)=\sum_{k=0}^{\infty}\ln\dfrac{(x+k+b-1)(k+1-b)}{(x+k)(k+1)}$, which is equivalent to the announced result.

REMARK 2.2.
(i) In the above proof we have interchanged heuristically the order of summation with differentiation and integration; still this interchange is in fact legitimate because the polynomial series (in x) involved, converges uniformly on any closed subinterval of $(0, +\infty)$, due to the very own construction of the $g_n$'s via the power series method.
In addition it is evident that $\Gamma(x+b)/\Gamma(x)$ is real analytic (in x) on $(0,\infty)$.
(ii) It is also evident that the infinite product representation of $f(1-b,b)$ should be considered equal to 1, a convention compatible with the statement of Proposition 2.1 and due to the obvious fact that $g_n=0$ for all n.

(iii) For notational convenience we will set the product of the first m terms in the infinite product representation of f(x,b) as $f_m(x,b)$ which we will call the <u>m-truncate</u> of the joint factor.

Using the functional equation of Proposition 2.1 for b=1/2 we obtain:

COROLLARY 2.2. (Legendre's formula revisited)
$$\Gamma(2x) = \frac{2^{2x-1}}{\pi} f(x,1/2)[\Gamma(x)]^2, \quad x>0 \qquad (2.8)$$

Euler's first integral ([4]) defines the Beta function which is well known to satisfy the functional equation B(x,y)=Γ(x)Γ(y)/Γ(x+y) for x>0, y>0; in the open strip (0,∞) x (0,1), by use of the joint factor f(x,y) and the infinite product representation of sin(πy), we can obtain immediately:

COROLLARY 2.3.
$$B(x,y) = (1/y) \prod_{k=1}^{\infty} [k(k-1+x+y)/(k+y)(k-1+x)]$$

REMARK 2.3. Corollary 2.2 provides us with a closed-type formula that not only recovers, but also simplifies many infinite product representations of $[\Gamma(x)]^2$, scattered throughout the pertinent literature ([1]). For example, it provides us with a representation of $[\Gamma(1/4)]^2$ "simpler" than the classical one (cf. [5]), which states (after squaring off) that $[\Gamma(1/4)]^2 = 4\pi \prod_{k=1}^{\infty} \frac{4k-1}{4k+1} \cdot \left(\frac{2k+1}{2k-1}\right)^{1/2}$. The aforementioned new representation reads $[\Gamma(1/4)]^2 = \pi(2\pi)^{1/2} \prod_{k=1}^{\infty} \frac{2k-1}{2k} \cdot \frac{4k-1}{4k-3}$, and it is a direct application of (2.8). The classical formula has a better "starting point" and, as simple algebra can demonstrate, it always "stays ahead" of its new counterpart. Still one gets compensated by the mere fact that the latter involves exclusively rational fractions and that for very large k's we achieve numerical approximations of the same order.

REMARK 2.4. Based on the functional equation Γ(1-x)= -xΓ(-x), which is valid for non integral x>0, we now see that Corollary 2.2 can easily be extended to negative arguments in the following form:

$$\Gamma(-2x) = -[\,2^{-2x}\text{cosec}(2\pi x)/xf(x,1/2)\,][\pi/\Gamma(x)]^2 \text{ for } 0<x<1/2.$$

REMARK 2.5. One could certainly had appealed to the Weierstrass canonical product representation of Γ(x+b) itself (cf. [3]) in order to achieve an ad hoc formula along the lines of Proposition 2.1; still the factorization in Proposition 2.1 and the discrete formal process that led to the latter, along with the emergence of the concept of the joint factor, seems to isolate a significant part of the "core" of this canonical product representation of Γ(x+b) without the direct involvment of non-elementary constants (like the Euler-Mascheroni constant).

## 3. A CLOSED-TYPE FORMULA FOR Γ(x), 1/Γ(x) AND ln Γ(x), x RATIONAL

For evident reasons, we will restrict ourselves to (0,1). We first prove a lemma which, by use of the joint factors introduced in Proposition 2.1, leads to the main result of this chapter:

LEMMA 3.1.
$$[\Gamma(1/p)]^p = \frac{(2\pi)^{p-1}}{p} \prod_{k=1}^{p-2} \frac{1}{f(1/p, k/p)}, \quad \text{for p any integer} \geq 3. \quad (3.1)$$

PROOF. By use of Gauss' multiplication formula and setting n=p and x=1/p in (1.4), we obtain in succession, via our functional equation of Proposition 2.1, that

$$\Gamma(1/p) = \Gamma(1/p)$$
$$\Gamma(2/p) = f(1/p, 1/p) \frac{\Gamma(1/p)}{\Gamma(1-1/p)}$$
$$\Gamma(3/p) = f(1/p, 2/p) \frac{\Gamma(1/p)}{\Gamma(1-2/p)}$$
$$\ldots\ldots$$
$$\Gamma((p-1)/p) = f(1/p, (p-2)/p) \frac{\Gamma(1/p)}{\Gamma(1-(p-2)/p)}.$$

Multiplying termwise the above equations and since: $\prod_{k=1}^{p-1} \Gamma(k/p) = \left[\frac{(2\pi)^{p-1}}{p}\right]^{1/2}$ we obtain the announced result.

The functional equation (1.2) and the above lemma gives now directly the following:

COROLLARY 3.1. $\Gamma(-1/p)^p = [-1/\sin(\pi/p)]^p [\pi/(2p)^{p-1}] \prod_{k=1}^{p-2} f(1/p, k/p)$.

We are ready now to prove the subsequent
THEOREM 3.1
$$\Gamma(q/p) = C_{p,q} \prod_{k=1}^{q-1} f(k/p, 1/p) \prod_{k=1}^{p-2} \left[\frac{1}{f(1/p, k/p)}\right]^{q/p} \quad (3.2)$$

where $C_{p,q} = \frac{(2\pi)[2\sin(\pi/p)]^{q-1}}{(2\pi p)^{q/p}}$ for p, q integers and $1 < q \leq p$, $p \geq 3$.

PROOF OF (3.2). Following the same idea as in Lemma 3.1, and by repeated use of Proposition 2.1, we have that:
$$\Gamma(1/p) = \Gamma(1/p)$$
$$\Gamma(2/p) = \Gamma(1/p) \frac{f(1/p, 1/p)}{\Gamma(1-1/p)}$$
$$\Gamma(3/p) = \Gamma(2/p) \frac{f(2/p, 1/p)}{\Gamma(1-1/p)}$$
$$\ldots\ldots\ldots$$
$$\Gamma(q/p) = \Gamma((q-1)/p) \frac{f((q-1)/p, 1/p)}{\Gamma(1-1/p)}$$
$$\Gamma(q/p) = \left[\frac{\sin(\pi/p)}{\pi}\right]^{q-1} [\Gamma(1/p)]^q \prod_{k=1}^{q-1} f(k/p, 1/p)$$

and thus, by use of (3.1), we arrive at the announced result.

One can now automatically obtain the following results:

COROLLARY 3.2.

(i) $\dfrac{1}{\Gamma(q/p)} = \dfrac{1}{C_{p,q}} \prod_{k=1}^{p-2} [f(1/p,k/p)]^{q/p} \prod_{k=1}^{q-1} \dfrac{1}{f(k/p,1/p)}$

(ii) $\ln\Gamma(q/p) = (\ln C_{p,q}) + \sum_{k=1}^{q-1} \mu_k(p) - (q/p) \sum_{k=1}^{p-2} v_k(p)$, where $\mu_k(p) = \ln f(k/p, 1/p)$, $1 \leq k \leq q-1$

and $v_k(p) = \ln f(1/p, k/p)$, $1 \leq k \leq p-2$.

REMARK 3.1
   (a) It is now evident that, via the classical continuity argument, (3.2) provides in principle the mechanism for the evaluation of Γ(x) at any number x in (0,1) and thus at any real x not in $Z_0^-$. But the actual numerical implementation of this remark would demand, aside to a sophisticated programming, the aid of numerical analysis methods and the pertinent approximation theory and thus lies beyond the scope of the present work.
   (b) Theorem 3.1 provides an alternative way to reproduce a rather big part of the explcit formulae presented in Section 2 of [9] (pp 269-270).

## 4. APPLICATIONS

In this chapter we demonstrate the method that allows us, by use of the algorithm (2.5) as modified in Remark 2.1, to produce interesting representations of non classical values of Γ(x) and other transcendental functions, as well as some presumably new inequalities and features concerning Γ(x) and its applications[*]:

APPLICATION 1.
It is rather evident that the joint factor f(1/p,k/p) has $f_m$(1/p,k/p) as a strict lower bound for all positive integers p, k and m such that p>2 and k≥p-2; in particular, for m=1, Lemma 3.1 and Corrolary 3.1 combined lead now to the following pair of inequalities:
(i) $\Gamma(1/p) < (2\pi/p)^{1-(1/p)} [(p-1)!]^{2/p}$

(ıı) $\Gamma(-1/p) > - [\pi/\sin(\pi/p)] (2\pi)^{1-(1/p)}$

Note that, numerically, (i) and (ii) provide rather crude bounds because of the particularly small value of m; the larger the value of m the more refined inequalities will be produced. We illustrate this point extensively in the reasoning of Applications 9 & 10.

APPLICATION 2.

For x>0, using Proposition 2.1 twice and the resulting product of the two joint factors as a single infinite product, we see that $x\,\Gamma(x) = \Gamma(x+1) = \Gamma(x+1/2)\,f(x+1/2,1/2)/\sqrt{\pi}$ = $\Gamma(x)f(x,1/2)f(x+1/2,1/2)/\pi$ which results to $\prod_{k=1}^{\infty}(2k/2k-1)^2[(k-1+x)/k+x] = \pi x$; the latter when combined with the infinite product representation of $(\sin\pi x)/\pi x$ produces, another representation of $\sin(\pi x)$, namely $\sin(\pi x) = \prod_{k=1}^{\infty}(2/2k-1)^2[(k-x)/(k-1+x)]$.

APPLICATION 3.
Setting x=b, 0 < b < 1, in Proposition 2.1 and using Corrolary 2.2 one gets immediately that $\sin(\pi b)f(b,b)/f(b,1/2) = 2^{2b-1}$; then employing the joint products involved, we obtain the subsequent formula that, modulo the formulae that one could randomly encounter when surfing on the web, we were not able to trace in the pertinent literature: $\prod_{n=1}^{\infty}(n-1/2)(n-1/2+b)/(n-b)(n-1+2b) = 2^{2b-1}/(\sin\pi b)$.

Note that in particular for b=1/4 and b=1/2 we set by convention the corresponding infinite product of 1's equal to 1.

APPLICATION 4.

By Proposition 2.1 for b=1/2 we obtain that $\Gamma\left(x+\dfrac{1}{2}\right) = \dfrac{f(x,1/2)}{\pi^{1/2}}\Gamma(x)$; we also have that $\Gamma(1-x) = \Gamma\left(\dfrac{1}{2}+\dfrac{1}{2}-x\right) = \dfrac{f(1/2-x,1/2)}{\pi^{1/2}}\Gamma\left(\dfrac{1}{2}-x\right)$ for 0<x<1/2. We conclude that

$\Gamma(1/2-x) = \dfrac{\Gamma(1-x)\pi^{1/2}}{f(1/2-x,1/2)}$ and thus $\Gamma(1/2+x)\Gamma(1/2-x) = \dfrac{\pi}{\cos(\pi x)} =$

$= \dfrac{f(x,1/2)}{f(1/2-x,1/2)}\Gamma(x)\Gamma(1-x) = \dfrac{f(x,1/2)}{f(1/2-x,1/2)}\dfrac{\pi}{\sin(\pi x)}$. We see now that we have

produced a seemingly new infinite product representation for $\tan(\pi x)$, namely

$\tan(\pi x) = \prod_{n=1}^{\infty}\dfrac{(2n-2+2x)(2n-2x)}{(2n-1)^2-(2x)^2}$. The corresponding classical infinite product

representation -directly from the classical math literature, or indirectly when one combines the formula stated for $\sin(\pi x)$ with the corresponding one for $\cos(\pi x)$, (e.g. 4.3.90 in [6])- produces the "more complicated" classical

$\tan(\pi x) = \pi x \prod_{n=1}^{\infty}[(2n-1)/n]^2\left[(n^2-x^2)/(2n-1)^2-(2x)^2\right]$.

APPLICATION 5.
Let us combine now the basic fact of the strict convexity of $\Gamma(x)$ over $(0,+\infty)$ with the factorization presented in Proposition 2.1 in order to produce a strict upper bound for $\Gamma(\alpha)$, 0<α<1, provided α stays away from zero. We know that the classical inequality $\Gamma(\alpha x+(1-\alpha)y) < \alpha\Gamma(x)+(1-\alpha)\Gamma(y)$, is true for x>0, y>0, x≠y; if we set x=1+y and y=1-α+ε with 0<ε<α<1, then using the fact that the joint factor f(1-α,ε) is strictly bounded from above by any of its m-truncation (compare to the similar reasoning of Application 2), we obtain:

$$\Gamma(\alpha) < [\sin(\pi\varepsilon) / \sin(\pi\alpha)][(1+\varepsilon\alpha^2)f_m(1-\alpha,\varepsilon) / \varepsilon] \text{ for all } m \in \mathbb{N}.$$

In particular, when we set m=1, we have the following inequality:

$\Gamma(\alpha) < [\sin(\pi\varepsilon) / \sin(\pi\alpha)][(1-\alpha)(1+\varepsilon\alpha-\alpha^2) / \varepsilon(1-\varepsilon)(1+\varepsilon-\alpha)]$; letting $\varepsilon=\alpha$ we have the indicative approximation $\Gamma(\alpha) < 1/\alpha$ on (0,1), which evidently performs quite satisfactory whwnever α is close to the endpoints. It is also clear that for fixed ε and α (ε≠α), the larger the m the better the m-ub for $\Gamma(\alpha)$.

APPLICATION 6.
It is immediate now that Cor.2.4 leads to the following inequalities:

(i) If 0<x,y<1 and m ∈ N, then $B(x,y) > (1/y) \prod_{k=1}^{m} [k(k-1+x+y)/(k-1+x)(k+y)]$

(ii) If x>1, 0<y<1 the inequality in (i) is reserved.

APPLICATION 7.
Let $\Omega_n$ be the volume of the unit ball in $\mathfrak{R}^n$, n ∈ N. Then since $\Omega_n = \pi^{n/2} / \Gamma(1+n/2)$ we immediately conclude that $\Omega_{n-1}/\Omega_n = f((n+1)/2,1/2) / \pi = n / [2f(n/2,1/2)]$. We can directly obtain, for all n, m, s ∈ N, the following inequalities:

$$(1/\pi) \prod_{k=1}^{s} (2k/2k-1)[(2k+n-1)/(2k+n)] < \Omega_{n-1}/\Omega_n \leq$$

$$(n/2) \prod_{k=1}^{m} [(2k-1)/2k][(2k+n-1)/(2k+n-2)].$$

Note that the right-hand side equality holds iff n=1; these s-l.bounds and m-u.bounds seem to be the simplest rational bounds (in s and m respectively) to hold for all n but they are less sharp for each n, whenever compared to the corresponding ones provided by Borgwardt in [5].

APPLICATION 8.
The Wallis's integral is defined as $I_\alpha = \int_0^{\pi/2} (\sin\theta)^\alpha d\theta = \int_0^{\pi/2} (\cos\theta)^\alpha d\theta$. For α>-1 we know that $I_\alpha = (1/2) B((\alpha+1)/2, 1/2)$ ([1]). Using the joint factor concept it is now immediate that if α>0 then $I_\alpha = (1/\alpha) f(\alpha/2, 1/2)$. Then using the m-truncation concept we can easily deduce the following inequalities:

(i) For 0<α<1 and any m ∈ N, $I_\alpha < [2^m (m!)(\alpha/2)_m] / [(2m-1)!!((\alpha+1)/2)_m]$, where $(.)_m$ denotes Pochammer's symbol ([6]).
(ii) For α>1 the inequality in (i) is reversed.

Note that even in the case of the crudest possible estimates, obtained for m=1, the above inequalities provide a better l.b. and u.b., respectively, of Wallis's integral in the corresponding interval than the ones provided by the classical trigonometric inequalities. For example, when 0<θ<π/2, we know that 2θ/π < sinθ < θ and thus we conclude that $\pi/2(\alpha+1) < I_\alpha < (\pi/2)^{\alpha+1} / (\alpha+1)$. On the other hand, for 0<α<1, (i) states that:

$$I_\alpha < 2/(\alpha+1) < (\pi/2)^{\alpha+1} /(\alpha+1),$$

while for α>1 (ii) states that:

$$I_\alpha > 2/(\alpha+1) > \pi/2(\alpha+1).$$

place accuracy; on the other hand the "old classical formula" for the same upper index produces exactly the same accuracy.

APPLICATION 9.

In ([2]) H. Alzer, among other sharp inequality conditions about the Gamma function, obtained the following improvement of Gamma inequalities discovered in 1997 by Anderson and Qiu:
Let $\alpha=1-\gamma$ and $\beta=(\pi^2-6\gamma)/12$, where $\gamma$ denotes Euler's constant. Then
(i) If $0<x<1$, $A(x) = x^{\alpha x-1} < \Gamma(x) < D(x) = x^{\beta(x-1)-\gamma}$.
(ii) If $1<x<\infty$, $D(x) < \Gamma(x) < E(x) = x^{x-1-\gamma}$.
Using the joint factor $f(x,1/2)$ and its m-truncation we see that
(a) If $0<x<1/2$ then:
$$\Gamma(x) > \left(\sqrt{\pi}\right)[\Gamma(x+1/2)]/f_m(x,1/2) > \left(\sqrt{\pi}\right)A(x+1/2)/f_m(x,1/2) = K_m(x).$$
(b) If $1/2<x<\infty$ then:
$$\Gamma(x) < \left(\sqrt{\pi}\right)[\Gamma(x+1/2)]/f_m(x,1/2) < \left(\sqrt{\pi}\right)E(x+1/2)/f_m(x,1/2) = L_m(x).$$

It emerges now as a rather natural question whether these l.b.'s (resp. u.b.'s) in (a) (resp. (b)) refine, for a particular choice of m, Alzer's Γ-bounds given in (i) and/or (ii). We can easily check that each term of the functional sequence $(K_m(x))$ when restricted to a suitable subinterval $I_m$ of $(0,1/2)$, refines $A(x)$; also each term of $(L_m(x))$ when restricted to a suitable subinterval $J_m$ of $(1/2,1)$, refines $D(x)$, while when restricted to a suitable subinterval $P_m$ of $(1, \infty)$, refines $E(x)$. It is evident by their construction and via a simple continuity argument that these intervals are nested and their lengths strictly increase to 1/2, 1/2 , and ∞, respectively, as m→ ∞.
We provide here three numerical examples that illuminate the above observation for the simplest in form l.b. and u.b. refinement of Alzer's bounds in predetermined subintervals. To that end we restrict
(a) the gamma l.b. $K_1(x) = \left(\sqrt{\pi}/2x\right)(x+0.5)^{\alpha(x+0.5)}$ to $I_1 = [0.241, 0.5)$,
(b) the gamma u.b. $L_1(x) = \left(\sqrt{\pi}/2x\right)(x+0.5)^{x+0.5-\gamma}$ to $J_1 = (0.5, 0.526]$,
and (c) the gamma u.b. $L_1(x)$ to $P_1 = [1.562, \infty)$.
(Here for simplicity we have settled for a third decimal place accuracy concerning the interval endpoints.)
After the evident substitutions, the refinement inequalities to be demonstrated are, respectively, equivalent to:
(d) $h(x) = (1+0.5/x)^x \sqrt{(x+0.5)} > \delta = \left(2/\sqrt{\pi}\right)^{1/\alpha}$ on $I_1$.
(e) $g(x) = (x+0.5)^{x+0.5-\gamma} / x^{\beta(x-1)+\alpha} < \varepsilon = 2/\sqrt{\pi}$ on $J_1$.
(f) $s(x) = (x+0.5)^{x+0.5-\gamma} / x^{x+\gamma} < \varepsilon$ on $P_1$.
Proof of (d): we can check that $h(0.240) < \delta < h(0.241)$; it suffices thus to show that $h(x)$ is strictly increasing on $(0,0.5)$. We see that $h(x)' = h(x)w(x)$, with $w(x)=\ln(1+0.5/x)-0.5/x+0.5/(x+0.5)$. But $w(x)' = (-0.5x^2+0.25x+0.125)/x^2(x+0.5)^2 > 0$ and $w(0^-) = 0$. Q.E.D.
Proof of (e): similarly we check first that $g(0.526) < \varepsilon < g(0.527)$ and show that $g(x)$ is strictly increasing of $(0.5,1)$: $g(x)' = g(x)[v(x)+r(x)]$ with $v(x) = [\ln(x+0.5) - \beta\ln x + 1 - \beta]$ and $r(x) = \gamma / (x+0.5) + (\beta-\alpha) / x$. It will suffice to show that min $v(x) >$ max $r(x)$. But min$[\ln(x+0.5) - \beta\ln x] = 0.37\ldots$, while max $r(x) = 3\gamma + 2\beta -1$ and $0.37 > 3(\beta+\gamma-1)$. Q.E.D.
Proof of (f): we check that $s(1.562) < \varepsilon < s(1.561)$ and that $s(x)$ is strictly decreasing on $(1, \infty)$ which is true since $s(x)' = s(x) [\ln(1+0.5/x) - \gamma / (x+0.5) - \gamma / x] < s(x) [(0.5 - \gamma) / (x + 0.5) - \gamma / x] < 0$.
Q.E.D

APPLICATION 10.

We conclude this chapter by examining a numerical impact of the joint factor approach in the revived research for improvement of Stirling's formula; we will make use of Nörlund's asymptotic representation of the Gamma function ([7]) and of a relatively recent result due to W. Schuster ([8]).
Stirling's formula states that $\Gamma(x) = (2\pi)^{1/2} x^{x-1/2} e^{-x} e^{\mu(x)}$, with $0 < \mu(x) < 1/12x$, while Nörlund's more effective asymptotic relation reads $\Gamma(x+1/2)=(2\pi)^{1/2} x^x e^{-x} e^{v(x)}$, where the bounds of $v(x)$ have been improved (see (17) in [7]) as follows:
$$(-1/24)[(1/x) + (1/120x^3)] \le v(x) \le (-1/24)[(1/x) - (1/8x^3)].$$
Combining now the joint factor $f(x,1/2)$ and Stirling's formula in Nörlund's relation, we easily obtain that $v(x) = \mu(x) + \ln[f(x,1/2) / \sqrt{(\pi x)}]$ and thus for $0<x<1/2$ we obtain

$v(x) < \mu(x) + \ln[f_m(x,1/2) / \sqrt{(\pi x)}]$ for all m, while for $x>1/2$ this inequality is reversed.
In order to avoid cumbersome expressions as bounds and at the same time illuminate our point, we will intentionally "sacrifice" a rather significant part of the numerical effectiveness of these inequalities by using again m=1. Then we can see directly that:
(i) If $0<x<1/2$ then $v(x) < (1/12x) + \ln\left[4\sqrt{x} / (1+2x)\sqrt{\pi}\right]$.
(ii) If $x>1/2$ then $v(x) > \ln\left[4\sqrt{x} / (1+2x)\sqrt{\pi}\right]$.
We are going at this point to demonstrate that the bound in (i), though the most crude among our m-type upper bounds, improves strictly Schuster's u.b. of $v(x)$, ultimately improving Nörlund's asymptotic representation of the Gamma function on (0,1/2):
We have to prove on (0,1/2) that $(1/8x) + \ln\left[4\sqrt{x} / (1+2x)\sqrt{\pi}\right] - \left(\frac{1}{192} x^3\right) < 0$;

since $\sqrt{x}/1+2x < \sqrt{2}/4$ it will suffice to show that $f(x) = 96 [\ln(2/\pi)]x^3 + 24x^2 - 1 < 0$.
We can easily check now, using the derivative test, that $\limsup f(x) = f(0.5^-) < 0$.
Q.E.D.
REMARK. By default, we can improve also the l.b. of $\mu(x)$ in a suitable subinterval of (0,1/2); specifically $\mu(x) > -\ln\left[4\sqrt{x} / (1+2x)\right]\sqrt{\pi} + 1/24\left[(1/x) + \left(1/120x^3\right)\right] > 0$, for $0.144 \le x < 0.5$, thus improving the $\Gamma(x)$-l.b. provided by Stirling's formula in the above interval. On the other hand, it is elementary to show that our m-type lower bounds in (ii) do not refine Schuster's l.b. of $v(x)$ for any $x>1/2$ and any m and the above considerations do not improve Stirling's u.b. $1/12x$ of $\mu(x)$.

## 5. THE DIGAMMA AND THE TRIGAMMA FUNCTION.

We conclude this work by considering two more general applications of the concept of the joint product, as well as of algorithm (2.5) concerning the psi (or digamma) function $\psi(z)$, and the *t*rigamma function $\psi'(z)$ respectively. The strategy to be presented can be repeated with the analogous results in case of any polygamma function restricted to (0,1) and and it can lead to various rather tedious series whose formulation entails theoretical aspects also.
It is an immediate consequence of the proof of Theorem 2.1 that:

$\psi(x+b)=\psi(x) - \sum_{n=1}^{\infty} g'_n(x,b)$ for $x>0$, $0 \le b < 1$ and $-\sum_{n=1}^{\infty} g'_n(x,b) = \frac{f'(x,b)}{f(x,b)}$ (where we restrict b to (0,1) to avoid indeterminancy for the case b=0). If we resort, after differentiating (2.5) w.r.t to x, to the resulting new recursive formula that involves the sequence $g'_n(x,b)=h_n(x,b)$, we will be able to obtain $\psi(x+b)$ and also, whenever necessary, the form of f '(x, b). This result besides being helpful in producing and/or reproducing identities as well as approximating $\psi(x)$ and $\psi'(x)$ for x>0, seems to bear an interest of its own. We present this proposed algorithmic strategy by providing in detail an illuminating abstract example concerning the case x=1-b, 0<b<1.

At first we deal with the function $\psi(x)$. Let us recall that $g_1(1-b,b)=g_2(1-b,b)=0$. For n≥2 algorithm (2.5) states that:

$(n+1)^2 g_{n+1}(x, b) =$

$= n[n+1-x-b(1-b-x)]g_n(x, b)+ \sum_{k=0}^{n-2}(k+1)g_{k+1}(x,b)[(n-k-1)g_{n-k-1}(x,b)-(n-k)g_{n-k}(x,b)]$

and thus $g_n(1-b, b)=0$ for all n. Preserving the x notation described in Remark 2.1 and after differentiating wrt x, we obtain at the point x=1-b

$$(n+1)^2 h_{n+1} = n(n+b)h_n \qquad (5.1)$$

for n≥2, where $h_1 = -b$ and $h_2 = -\frac{1}{4}b(1+b)$; note that we opressed the dependence of $h_n$ from b for reasons of simplicity. Multiplying (5.1) for n≥2 and after telescoping we arrive at

$$h_n = -\frac{1}{n(n!)}\frac{\Gamma(n+b)}{\Gamma(b)} \qquad (5.2)$$

for n≥3. Observe also that for n=1,2 we also obtain the correct values for $h_n$ and thus (5.2) holds for all n. Since $\psi(x+b) = \psi(x) - \sum_{n=1}^{\infty} h_n(x,b)$ for all x>0, 0<b<1, we derive that $\psi(1) = -\gamma = \psi(1-b) + \sum_{n=1}^{\infty} \frac{1}{n(n!)}\frac{\Gamma(n+b)}{\Gamma(b)}$. We conclude that for 0 < t < 1

$$\psi(t) = -\gamma - \frac{1}{\Gamma(1-t)}\sum_{n=1}^{\infty} \frac{\Gamma(n+1-t)}{n(n!)}. \qquad (5.3)$$

In terms of joint factors we have that $\Gamma(n+1-t) = \frac{(n-1)!}{\Gamma(t)}\cdot f(n,1-t)$ and thus (5.3) leads to:

$$\psi(t) = -\gamma - \frac{\sin(\pi t)}{\pi}\sum_{n=1}^{\infty} \frac{f(n,1-t)}{n^2} \qquad (5.4)$$

REMARK 5.1. A rather interesting improvement of the precision that can be obtained from (5.4) is feasible by use of the evident fact that $f(n,1-t) = o(\Gamma(t)n^{1-t})$:

Suppose that $n_0$ is very large in the identity $\sum_{n=1}^{\infty}\frac{f(n,1-t)}{n^2} = \sum_{n=1}^{n_0}\frac{f(n,1-t)}{n^2} + \sum_{n=n_0+1}^{\infty}\frac{f(n,1-t)}{n^2}$.

Thus we see that $\sum_{n=1}^{\infty}\frac{f(n,1-t)}{n^2} \approx \sum_{n=1}^{n_0}\frac{f(n,1-t)}{n^2} + \Gamma(t)\sum_{n=n_0+1}^{\infty}\frac{1}{n^{1+t}}$ and so we conclude that:

$$\psi(t) \approx -\gamma - \frac{\sin(\pi t)}{\pi}\sum_{n=1}^{n_0}\frac{f(n,1-t)}{n^2} - \frac{1}{\Gamma(1-t)}\sum_{n=n_0+1}^{\infty}\frac{1}{n^{1+t}} \qquad (5.5)$$

REMARK 5.2. It is rather evident that for $n_0$ very large (5.5) furnishes one more "interconnection" between the psi function and the Riemann's zeta function when the former is restricted to (0,1), with the aid of the joint factors namely:

$$\psi(t) \approx -\gamma + \left[\frac{1}{\Gamma(1-t)} - \frac{\sin(\pi t)}{\pi}\right] \sum_{n=1}^{n_0} \frac{f(n,1-t)}{n^2} - \frac{\zeta(1+t)}{\Gamma(1-t)}.$$

In order to obtain a similar expression for the trigamma function, let us now differentiate (5.4) wrt t. We then have that:

$$\psi'(t) = -\cos(\pi t) \sum_{n=1}^{\infty} \frac{f(n,1-t)}{n^2} - \frac{\sin(\pi t)}{\pi} \sum_{n=1}^{\infty} \frac{f'(n,1-t)}{n^2}.$$

Since we are dealing once more with the derivative of an infinite product of functions (in t) we eleviate this difficulty through the classical logarithmic differentiation, namely $f'(n,1-t) = -\sum_{k=1}^{\infty} \frac{f(n,1-t)}{k-1+t}$. By substitution in the above formula of $\psi'(t)$ we arrive at:

$$\psi'(t) = \frac{\sin(\pi t)}{\pi} \sum_{n=1}^{\infty} \frac{f(n,1-t)\sum_{k=1}^{n}\frac{1}{k-t}}{n^2} \qquad (5.6)$$

Finally let us construct another approximation formula that will improve the calculating efficiency of (5.6) working along the same lines as with (5.5); we can immediately conclude that, for $n_0$ very large:

$$\psi'(t) \approx \frac{\sin(\pi t)}{\pi} \sum_{n=1}^{n_0} \frac{f(n,1-t)\sum_{k=1}^{n}\frac{1}{k-t}}{n^2} + \frac{1}{\Gamma(1-t)} \sum_{n=n_0+1}^{\infty} \frac{\sum_{k=1}^{n}\frac{1}{k-t}}{n^{1+t}} \qquad (5.7)$$

### REFERENCES
1. M. Abramowitz and I.A. Stegun, Handbook of Mathematical Functions, Dover, New York, 1972.
2. H. Alzer, Inequalities for the Gamma Function, Proc. A.M.S (2000), v.128, 141-147.
3. G. Andrew's, R. Askey and R. Roy, Special Functions, Cambridge Univ. Press, 1999.
4. E. Artin, The Gamma Function, Holt, Rinehart and Winston (eds), 1964.
5. K.H. Borgwardt, The Simplex Method, Springer –Verlag, Berlin,1987.
6. I.S. Gradshteyn and I.M. Ryzhik, Table of Integrals, Series and Products, (Corr.& Enlarg. by Alan Jeffrey), Academic Press, 1980.
7. D.Karayannakis and I.S.Xezonakis, A new algorithm for the Gamma function at rational points and its numerical ramifications (submitted manuscript)
8. W. Schuster, Improving Stirling's formula, Arch. Math, v.77 (2001),170 -176.
9. R.Vidunas, Expressions For Values of the Gamma Function, Kyushu Journal of Math.,Vol.59, 2005, pp 267-283
10. E.T. Whittaker and G.N. Watson, Modern Analysis, 4th ed., Cambridge Univ. Press, 1934.